\documentclass[secthm,seceqn,amsthm,ussrhead,10pt]{amsart}
\usepackage{amsmath,latexsym}
\usepackage[english]{babel}
\usepackage[psamsfonts]{amssymb}
\usepackage{times}
\usepackage{cite}
\usepackage{pdflscape} 
\usepackage{ulem}
\usepackage[mathcal]{euscript}
\usepackage{tikz}
\usepackage{hyperref}
\usepackage{cancel}
\usetikzlibrary{arrows}

\newtheorem{Th}{Theorem}
\newtheorem{Lem}[Th]{Lemma}

\newtheorem{corollary}[Th]{Corollary}

\newenvironment{Proof}[1][Proof.]{\begin{trivlist}
\item[\hskip \labelsep {\bfseries #1}]}{\flushright
$\Box$\end{trivlist}}

\begin{document}
	\sloppy
	
{\huge Degenerations of binary Lie  and nilpotent Malcev  algebras
\footnote{The authors were supported
by RFBR 	16-31-00004,
FAPESP  14/24519-8,
FAPESP 14/19521-3
and by R \& D 6.38.191.2014 of Saint-Petersburg State University, "Structure theory, classification, geometry, K-theory and arithmetics of algebraic groups and related structures".}
}

\medskip

\medskip

\medskip

\medskip
\textbf{Ivan Kaygorodov$^{a}$, Yury Popov$^{b}$, Yury Volkov$^{c,d}$}
\medskip

{\tiny
$^{a}$ Universidade Federal do ABC, CMCC, Santo Andr\'{e}, Brazil.

$^{b}$ Novosibirsk state university, Novosibirsk, Russia.

$^{c}$ Universidade de S\~{a}o Paulo, IME, Sao Paulo, Brazil.

$^{d}$ Saint Petersburg state university, Saint Petersburg, Russia.
\smallskip

    E-mail addresses:

    Ivan Kaygorodov (kaygorodov.ivan@gmail.com),
    
    Yury Popov (yuri.ppv@gmail.com),
    
    Yury Volkov (wolf86\_666@list.ru).

}

       \vspace{0.3cm}

{\bf Abstract.} 
We describe degenerations of four-dimensional binary Lie  algebras, and five- and six-dimensional nilpotent Malcev algebras  over $\mathbb{C}.$  
In particular, we describe all irreducible components of these varieties.

{\bf Keywords:} Malcev algebra, binary Lie algebra, nilpotent algebra, degeneration, rigid algebra
       \vspace{0.3cm}

       \vspace{0.3cm}

\section{Introduction}

       \vspace{0.3cm}

Degenerations of algebras is an interesting subject, which was studied in various papers (see, for example, \cite{B99,B05,G75,M79,M80,AOR05,CKLO13,R06,KE14,KE16,BC99,S90,GRH,GRH2,BB09,BB14,B03,gorb93}).
In particular, there are many results concerning degenerations of algebras of low dimensions from some variety defined by a set of identities.
One of important problems in this direction is the description of so-called rigid algebras. These algebras are of big interest, since the closures of their orbits under the action of generalized linear group form irreducible components of a variety under consideration
(with respect to Zariski topology). For example, the problem of finding rigid algebras was solved for
low-dimensional associative (see \cite{G75,M79,M80}), Leibniz (see \cite{AOR05,CKLO13,R06}), and Jordan (see \cite{KE14,KE16}) algebras.
There are significantly less works where the full information about degenerations was found for some variety of algebras.
This problem was solved for four-dimensional Lie algebras in \cite{BC99}, for nilpotent five- and six-dimensional Lie algebras in \cite{S90,GRH}, for  two-dimensional pre-Lie algebras in \cite{BB09}, and for three-dimensional Novikov algebras in \cite{BB14}.

The notions of Malcev and binary Lie (BL for short) algebras were introduced by Malcev in \cite{M55}.
The structure theory and some properties of Malcev algebras were studied by Kuzmin and other authors
(see, for example, \cite{kuzma2, sagle61, kuz71, kay14, K70,kp16,poji01,fil2}).
Note that any Lie algebra is a Malcev algebra and any Malcev algebra is a BL algebra.
Note also that any alternative algebra can be turned to a Malcev algebra by defining a new multiplication $[,]$ by  $[x,y]=xy-yx$.
Any Malcev algebra is a tangent algebra of a suitable locally analytic Moufang loop (see \cite{kuz71}). 

In this paper we give the full information about degenerations of BL algebras of dimension 4 and nilpotent Malcev algebras of dimensions $5$ and $6.$
More precisely, we construct a graph of primary degenerations. The vertices of this graph are isomorphism classes of algebras from the variety under consideration.
An algebra $A$ degenerates to an algebra $B$ iff there is a path from the vertex corresponding to $A$ to the vertex corresponding to $B$.
Thus, we obtain a generalization of analogous results of \cite{S90,BC99,GRH} for Lie algebras. Also we describe rigid algebras and irreducible components for these varieties of algebras.


\section{Definitons and notation}

All spaces in this paper are considered over $\mathbb{C}$, and we write simply $dim$, $Hom$ and $\otimes$ instead of $dim_{\mathbb{C}}$, $Hom_{\mathbb{C}}$ and $\otimes_{\mathbb{C}}$. An algebra $A$ is a set with a structure of vector space and a binary operation that induces a bilinear map from $A\times A$ to $A$.

Given an $n$-dimensional vector space $V$, the set $Hom(V \otimes V,V) = V^* \otimes V^* \otimes V$ 
is a vector space of dimension $n^3$. This space has the structure of the affine variety $\mathbb{C}^{n^3}.$ Indeed, let us fix a basis $e_1,\dots,e_n$ of $V$. Then any $\mu\in Hom(V \otimes V,V)$ is defined by structure constants $c_{i,j}^k\in\mathbb{C}$ such that
$\mu(e_i\otimes e_j)=\sum\limits_{k=1}^nc_{i,j}^ke_k$. A subset of $Hom(V \otimes V,V)$ is called closed if it can be defined by a set of polynomial equations in variables $c_{i,j}^k$ ($1\le i,j,k\le n$).

Let $T$ be a set of polynomial identities.
All algebra structures on $V$ satisfying polynomial identities from $T$ form a Zariski-closed affine subset of the variety $Hom(V \otimes V,V)$. We denote this subset by $\mathbb{L}(T)$.
The general linear group $GL_n(\mathbb{C})$ operates on $\mathbb{L}(T)$ by conjugation:
$$ (g * \mu )(x\otimes y) = g( \mu(g^{-1}(x)\otimes g^{-1}(y)))$$ 
for $x,y\in V$, $\mu\in \mathbb{L}(T)\subset Hom(V \otimes V,V)$ and $g\in GL_n(\mathbb{C})$.
Thus, $\mathbb{L}(T)$ is decomposed into $GL_n(\mathbb{C})$-orbits that correspond to the isomorphism classes of algebras. 
Let $O(\mu)$ denote the orbit of $\mu\in\mathbb{L}(T)$ under the action of $GL_n(\mathbb{C})$. Correspondingly, $\overline{O(\mu)}$ is the Zariski closure of $O(\mu)$.

Let $A$ and $B$ be two $n$-dimensional algebras satisfying identities from $T$. Let $\mu$ and $\lambda$ from $\mathbb{L}(T)$ represent $A$ and $B$ respectively.
We say that $A$ degenerates to $B$ and write $A\to B$ if $\lambda\in\overline{O(\mu)}$.
Note that in this case we have $\overline{O(\lambda)}\subset\overline{O(\mu)}$. Hence, the definition of degeneration does not depend on the choice of $\mu$ and $\lambda$. We write $A\not\to B$ if $\lambda\not\in\overline{O(\mu)}$.

Let $A$ be represented by the structure $\mu\in\mathbb{L}(T)$. The algebra $A$ is called {\it rigid} in $\mathbb{L}(T)$ if $O(\mu)$ is an open subsetset of $\mathbb{L}(T)$.
 Recall that a subset of a variety is called irreducible if it can't be represented as a union of two non-trivial closed subsets. A maximal irreducible closed subset of a variety is called an {\it irreducible component}.
 In particular, $A$ is rigid in $\mathbb{L}(T)$ iff $\overline{O(\mu)}$ is an irreducible component of $\mathbb{L}(T)$. Let $Irr(\mathbb{L}(T))$ and $Rig(\mathbb{L}(T))$ denote  the set of irreducible components of $\mathbb{L}(T)$ and the set of rigid algebras in $\mathbb{L}(T)$ respectively.  It is well known that any affine variety can be represented as a finite union of its irreducible components in a unique way.
 
Let $A$ be an algebra. For $x,y,z \in A$ we define their Jacobian $J(x,y,z)$ by the equality $J(x,y,z)=(xy)z+(yz)x+(zx)y.$
The algebra $A$ is called a Malcev algebra if it satisfies the identities
$$xy=-yx, \ J(x,y,xz)=J(x,y,z)x.$$
The algebra $A$ is called a binary Lie (BL for short) algebra if all its $2$-generated subalgebras are Lie algebras.
It was shown by Gainov in \cite{gai57} that $A$ is a BL algebra iff it satisfies the identities
$$xy=-yx, \ J(x,y,xy)=0.$$
It is easy to see that any Lie algebra is a Malcev algebra and any Malcev algebra is a BL algebra. It was shown in \cite{G63}
that any three-dimensional BL algebra is a Lie algebra.
The classification of four-dimensional BL algebras was obtained in \cite{G63,K98}.
The classification of five-dimensional and nilpotent six-dimensional non-Lie Malcev algebras is given in \cite{K70}.

Let $BL_n$, $Mal_n$ and $Lie_n$ denote the varieties of $n$-dimensional BL, Malcev and Lie algebras respectively, and $NBL_n$, $NMal_n$ and $NLie_n$ denote their subvarieties formed by nilpotent algebras. 

Define the sets $A^l$ by the equalities $A^1=A$ and $A^l=A^{l-1}A$ ($l>1$). Also define the central series $Z_l(A)$ ($l>0$) of $A$ in the following way.
We define $Z_1(A)=Z(A)$ as the center of $A$, and, for $l>1,$  $Z_l(A)$ is the full inverse image of $Z(A/Z_{l-1}(A))$ under the canonical projection from $A$ to $A/Z_{l-1}(A)$.

We collect all the information that we need about the algebras under consideration in Tables 4.1, 5.1 and 6.1. In these tables in the first column we write the names of the algebras. In the second column we give the multiplication tables in some fixed basis $e_1,\dots,e_n$ of $V$. All products of basis elements, which are not described in the table, are zero or can be deduced from one of the described products and the anticommutativity identity. In the third column we give the dimensions of algebras of derivations. In the columns named $Z(A)$, $A^2$ and $A^3$ we give the dimensions of the corresponding spaces. In the column named $Z_l(A)$ we give the dimensions of the members of central series of $A$.
Also, in the last column of Table 4.1 we have "Lie" for Lie algebras, "Malcev" for Malcev non-Lie algebras and "BL" for BL non-Malcev algebras.

The names of four-dimensional Lie algebras are from \cite{B03}. The classification of BL non-Lie algebras is taken from \cite{K98}.
One of them is called $g_3(\beta)$ here, since for $\beta=2$ we obtain the algebra $g_3$ in the notation of \cite{B03}.
The remaining BL non-Lie algebra is called $g_6$.

The names of five- and six-dimensional nilpotent Lie algebras are taken from \cite{S90}. The classification of Malcev non-Lie algebras of corresponding dimensions is deduced from \cite{K70}. We give names containing the letter "$M$" to these algebras. So in our notation a five- or six-dimensional algebra is Malcev and non-Lie iff it contains a letter "$M$" in its name, except the algebras $M_6^1$ and $M_2^1$ that correspond to the algebras $g_{6,4}$ and $g_{6,12}$ respectively in the notation of \cite{S90}.

\section{Methods} 

In the present work we use the methods that were applied for Lie algebras in \cite{BC99,GRH,GRH2,S90}.
First of all, it is well known that if $A\to B$ and $A\not\cong B$, then $dim\,Der(A)<dim\,Der(B)$, where $Der(A)$ is the algebra of derivations of $A$. We have computed the dimensions of algebras of derivations and have checked the assertion $A\to B$ only for such $A$ and $B$ that $dim\,Der(A)<dim\,Der(B)$. Secondly, it is well known that if $A\to C$ and $C\to B$, then $A\to B$. If there is no $C$ such that $A\to C$ and $C\to B$, then the assertion $A\to B$ is called a {\it primary degeneration}. If $dim\,Der(A)<dim\,Der(B)$ and there are no $C$ and $D$ such that $C\to A$, $B\to D$ and $C\not\to D$, then the assertion $A \not\to B$ is called a {\it primary non-degeneration}. It is enough to prove only primary degenerations and non-degenerations to describe all degenerations in the variety under consideration. It is easy to see that any algebra degenerates to the algebra with zero multiplication.

Degenerations of four-dimensional and nilpotent five- and six-dimensional Lie algebras were described in \cite{BC99,S90,GRH}.
Since the set $\mathbb{L}(T)$ is closed for any $T$, a Lie algebra can't degenerate to a non-Lie algebra.
So when we want to add Malcev or BL algebras to Lie algebras we don't have to check the degenerations from Lie algebras to any of the added algebras.

To prove the primary degenerations we construct the families of matrices parametrized by $t$. Namely, let $A$ and $B$ be two algebras represented by the structures $\mu$ and $\lambda$ from $\mathbb{L}(T)$ respectively. Let $e_1,\dots, e_n$ be a basis of $V$, for which  $\lambda$ is defined by structure constants $c_{i,j}^k$ ($1\le i,j,k\le n$). If there exist $a_i^j(t)\in\mathbb{C}$ ($1\le i,j\le n$, $t\in\mathbb{C}^*$) such that $E_i^t=\sum\limits_{j=1}^na_i^j(t)e_j$ is a basis of $V$ for $t\in\mathbb{C}^*$ and the structure constants of $\mu$ in the basis $E_1^t,\dots, E_n^t$ are such polynomials $c_{i,j}^k(t)\in\mathbb{C}[t]$ that $c_{i,j}^k(0)=c_{i,j}^k$, then $A\to B$. In this case  $E_1^t,\dots, E_n^t$ is called a {\it parametrized basis} for $A\to B$.

Tables 4.2 and 6.2 give parametrized bases for primary degenerations between four-dimensional BL algebras and six-dimensional nilpotent Malcev algebras respectively. These tables include all primary degenerations of the form $A\to B$, where $A$ is a non-Lie algebra.

We now describe the methods for proving primary non-degenerations. The main tool for this is the following lemma.

\begin{Lem}[\cite{B99,GRH}]\label{main}
Let $\mathcal{B}$ be a Borel subgroup of $GL_n(\mathbb{C})$ and $\mathcal{R}\subset \mathbb{L}(T)$ be a $\mathcal{B}$-stable closed subset.
If $A \to B$ and $A$ can be represented by a structure $\mu\in\mathcal{R}$, then there is a structure $\lambda\in \mathcal{R}$ representing $B$.
\end{Lem}

Since any Borel subgroup of $GL_n(\mathbb{C})$ is conjugate to the subgroup of upper triangular matrices, Lemma \ref{main} can be applied in the following way. Let $A$ and $B$ be two algebras. Let $\mu,\lambda$ be some structures in $\mathbb{L}(T)$ representing $A$ and $B$ respectively.
Suppose that there is a set of equations $Q$ in variables $x_{i,j}^k$ ($1\le i,j,k\le n$) such that if $x_{i,j}^k=c_{i,j}^k$ ($1\le i,j,k\le n$) is a solution of all equations from $Q$, then $x_{i,j}^k=\tilde c_{i,j}^k$ ($1\le i,j,k\le n$) is a solution for all equations from $Q$ too in the following cases:
\begin{enumerate}
    \item if $\tilde c_{i,j}^k=\frac{\alpha_i\alpha_j}{\alpha_k}c_{i,j}^k$ for some $\alpha_i\in\mathbb{C}^*$ ($1\le i\le n$);
    \item if there are some numbers $1\le u<v\le n$ and some $\alpha\in\mathbb{C}$ such that
    $$
    \tilde c_{i,j}^k=\begin{cases}
    c_{i,j}^k,&\mbox{ if $i,j\not=u$ and $k\not=v$},\\
    c_{u,j}^k+\alpha c_{v,j}^k,&\mbox{ if $i=u$, $j\not=u$ and $k\not=v$},\\
    c_{i,u}^k+\alpha c_{i,v}^k,&\mbox{ if $i\not=u$, $j=u$ and $k\not=v$},\\
    c_{i,j}^v-\alpha c_{i,j}^u,&\mbox{ if $i,j\not=u$ and $k=v$},\\
    c_{u,u}^k+\alpha (c_{v,u}^k+c_{u,v}^k)+\alpha^2c_{v,v}^k,&\mbox{ if $i=j=u$ and $k\not=v$},\\
    c_{u,j}^v+\alpha (c_{v,j}^v-c_{u,j}^u)-\alpha^2c_{v,j}^u,&\mbox{ if $i=u$, $j\not=u$ and $k=v$},\\
    c_{i,u}^v+\alpha (c_{i,v}^v-c_{i,u}^u)-\alpha^2c_{i,v}^u,&\mbox{ if $i\not=u$, $j=u$ and $k=v$},\\
    c_{u,u}^v+\alpha (c_{v,u}^v+c_{u,v}^v-c_{u,u}^u)+\alpha^2(c_{v,v}^v-c_{v,u}^u-c_{u,v}^u)-\alpha^3c_{v,v}^u,&\mbox{ if $i=j=u$ and $k=v$}.
    \end{cases}
    $$
\end{enumerate}
Assume that there is a basis $f_1,\dots,f_n$ of $V$ such that the structure constants of $\mu$ in this basis form a solution for all equations from $Q$, but there is no basis $\tilde f_1,\dots,\tilde f_n$ of $V$ such that the structure constants of $\lambda$ in it form a solution for all equations from $Q$. Then $A\not\to B$.

We will often use two particular cases of Lemma \ref{main}. Firstly, if $dim\,A^l<dim\,B^l$ for some $l>0$, then $A\not\to B$. Secondly,  if $dim\,Z_l(A)>dim\,Z_l(B)$ for some $l>0$, then $A\not\to B$. 
In the cases where these two criterions can't be applied, we define $\mathcal{R}$ by some conditions, which can be expressed in terms of a set of equations $Q$ satisfying the property described above, and give a basis for $V$, in which the structure constants of $\mu$ satisfy all equations from $Q$. We omit everywhere the verification of the fact that $Q$ satisfies the required conditions and the verification of the fact that structure constants of $\lambda$ in any basis do not satisfy some equation from $Q$. These verifications can be done by direct calculations.

Another argument for the non-degeneration that we use is the so-called $(i,j)$-invariant. Given $i,j>0$, we call $c_{i,j}\in\mathbb{C}$ an {\it $(i,j)$-invariant} for the algebra $A$ if
$$tr(ad \ x)^i \cdot tr(ad \ y)^j=c_{i,j}{tr((ad \ x)^i \circ  (ad \ y)^j)}$$
for all $x,y\in A$. If $c_{i,j}$ is an $(i,j)$-invariant for $A$, but at the same time it is not an $(i,j)$-invariant for $B$, then
$A\not\to B$.

We give the proof of primary non-degenerations in Tables 4.3 and 6.3, where for each primary non-degeneration we give one of the arguments mentioned above.

If the number of orbits under the action of $GL_n(\mathbb{C})$ on the variety $\mathbb{L}(T)$ is finite, then the graph of primary degenerations gives the whole picture. In particular, the description of rigid algebras and irreducible components can be easily obtained.
But in this work in some cases the situation is not so good. Then we have to be able to verify a little more complicated assertions. Let $A_*=\{A_{\alpha}\}_{\alpha\in I}$ be a set of algebras and $B$ be some other algebra. Suppose that $A_{\alpha}$ is represented by the structure $\mu_{\alpha}$ ($\alpha\in I$) and $B$ is represented by the structure $\lambda$. Then $A_*\to B$ means $\lambda\in\overline{\bigcup\limits_{\alpha\in I}O(\mu_{\alpha})}$, and $A_*\not\to B$ means $\lambda\not\in\overline{\bigcup\limits_{\alpha\in I}O(\mu_{\alpha})}$.

Let $A_*$, $B$, $\mu_{\alpha}$ ($\alpha\in I$) and $\lambda$ be as above. To prove that $A_*\to B$ we have to construct a family of pairs $(f(t), g(t))$ parametrized by $t$, where $f(t)\in I$ and $g(t)\in GL_n(\mathbb{C})$. Namely, let $e_1,\dots, e_n$ be a basis of $V$, for which $\lambda$ is defined by structure constants $c_{i,j}^k$ ($1\le i,j,k\le n$). If we construct $a_i^j(t)\in\mathbb{C}$ ($1\le i,j\le n$, $t\in\mathbb{C}^*$) and $f: \mathbb{C}^* \to I$ such that $E_i^t=\sum\limits_{j=1}^na_i^j(t)e_j$ is a basis of $V$ for $t\in\mathbb{C}^*$ and the structure constants of $\mu_{f(t)}$ in the basis $E_1^t,\dots, E_n^t$ are such polynomials $c_{i,j}^k(t)\in\mathbb{C}[t]$ that $c_{i,j}^k(0)=c_{i,j}^k$, then $A_*\to B$. In this case  $E_1^t,\dots, E_n^t$ and $f(t)$ are called a parametrized basis and a {\it parametrized index} for $A_*\to B$ respectively.

We now explain how to prove that $A_*\not\to B$. First of all, if $dim\,Der(A_{\alpha})>dim\,Der(B)$ for all $\alpha\in I$, then $A_*\not\to B$. One can use also the following generalization of Lemma \ref{main}, whose proof is the same as the proof of Lemma \ref{main}.

\begin{Lem}\label{gmain}
Let $\mathcal{B}$ be a Borel subgroup of $GL_n(\mathbb{C})$ and $\mathcal{R}\subset \mathbb{L}(T)$ be a $\mathcal{B}$-stable closed subset.
If $A_* \to B$, and for any $\alpha\in I$ the algebra $A_{\alpha}$ can be represented by a structure $\mu_{\alpha}\in\mathcal{R}$, then there is a structure $\lambda\in \mathcal{R}$ representing $B$.
\end{Lem}

\section{Binary Lie algebras of dimension 4}

The following table contains the classification and some invariants of four-dimensional BL algebras. It collects results from \cite{B03,K98}.

\begin{center}
\begin{equation*}
\begin{array}{|c|c|c|c|c|c|} 

\hline  A                       &\mbox{multiplication table}& Der(A)& Z(A) & A^2 & \mbox{type}  \\
\hline 
\hline  n_3 \oplus \mathbb{C}   & \begin{array} {c} e_1e_2=e_3\end{array} & 10      & 2    & 1   & Lie     \\
\hline  n_4                     & \begin{array} {c} e_1e_2=e_3, e_1e_3=e_4\end{array} &  7     & 1    & 2   & Lie \\
\hline  r_2 \oplus \mathbb{C}^2 & \begin{array} {c} e_1e_2=e_2\end{array}  & 8     & 2    & 1   & Lie     \\
\hline  r_2 \oplus r_2          & \begin{array} {c} e_1e_2=e_2, e_3e_4=e_4\end{array} & 4     & 0    & 2   & Lie     \\
\hline  sl_2 \oplus \mathbb{C}  & \begin{array}  {c} e_1e_2=e_2,  e_1e_3=-e_3, e_2e_3=e_1\end{array} & 4     & 1    & 3   & Lie     \\
\hline  g_1                     & \begin{array}  {c} e_1e_2=e_2,  e_1e_3=e_3,  e_1e_4=e_4\end{array} & 12    & 0    & 3   & Lie     \\
\hline  g_2(\beta)             & \begin{array}  {c} e_1e_2=e_2, e_1e_3=e_3, e_1e_4=e_3+\beta e_4\end{array}  & 8     & 0    & \begin{array}  {c}  3,\beta\neq0; \\ 2, \beta=0 \end{array}     & Lie     \\
\hline g_3(\beta)               & 
\begin{array}  {c} e_1e_2=e_2, e_1e_3=e_3, \\ e_1e_4=\beta e_4, e_2e_3= e_4\end{array}  &7     & 0    & 3   &   \begin{array}{c} \mbox{ $Lie$, for } \beta=2, \\
\mbox{ $Malcev$, for } \beta=-1 \\ \mbox{ $BL$, for } \beta\neq -1,2 \end{array}   \\

\hline  g_4(\alpha,\beta)       & \begin{array}  {c}  e_1e_2=e_2, e_1e_3=e_2+\alpha e_3, \\ e_1e_4=e_3+\beta e_4\end{array} & 6     & 0    & \begin{array}  {c} 3, \alpha\neq0\neq \beta; \\ 2, \alpha\beta=0 \end{array}      & Lie     \\
\hline  g_5(\alpha)             & \begin{array}  {c} e_1e_2=e_2,  e_1e_3=e_2+\alpha e_3, \\ e_1e_4=(\alpha+1)e_4,  e_2e_3=e_4\end{array} & 5     & 0    & \begin{array}  {c} 3, \alpha\neq0 \\ 2, \alpha=0 \end{array}    & Lie  \\

\hline g_6                        & \begin{array} {c} e_1e_2=e_3, e_3e_4=e_3\end{array} & 7     & 0    & 1   & BL     \\
\hline
 \end{array}
\end{equation*}
\end{center}
\begin{center} Table 4.1. Binary Lie algebras of dimension 4.
\end{center}

The algebra $g_4(\alpha_1,\beta_1)$ is isomorphic to $g_4(\alpha_2,\beta_2)$ iff
    the proportions $1:\alpha_1:\beta_1$ and $1:\alpha_2:\beta_2$ coincide after some permutation.
    The algebra $g_5(\alpha)$ is isomorphic to $g_5(\beta)$ iff $\alpha\beta=1$ or $\alpha=\beta.$ Apart from these two exceptions, any two algebras with different names from Table 4.1 are not isomorphic.

\begin{Th}\label{first}
The graph of primary degenerations for binary Lie  algebras of dimension 4 has the following form:
\end{Th}
\begin{center}

\begin{tikzpicture}[->,>=stealth',shorten >=0.05cm,auto,node distance=1cm,
                    thick,main node/.style={rectangle,draw,fill=gray!10,rounded corners=1.5ex,font=\sffamily \scriptsize \bfseries },rigid node/.style={rectangle,draw,fill=black!20,rounded corners=1.5ex,font=\sffamily \scriptsize \bfseries },style={draw,font=\sffamily \scriptsize \bfseries }]

\node (104)   {$4$};
\node (105) [below  of=104]      {$5$};
\node (106) [ below         of=105]      {$6$};
\node (107) [ below         of=106]      {$7$};
\node (108) [ below          of=107]      {$8$};
\node (110) [ below         of=108]      {$10$};
\node (112) [ below         of=110]      {$12$};
\node (116) [ below         of=112]      {$16$};

\node (1041)[right of=104]{};	
\node (1040)[right of=1041]{};	
	\node[rigid node] (1)  [right of =1040]                        {$r_2 \oplus r_2$ };
\node (1042)[right of=1]{};
\node (1043)[right of=1042]{};
\node (1044)[right of=1043]{};
\node (1045)[right of=1044]{};
\node (1046)[right of=1045]{};
\node (1047)[right of=1046]{};
	\node[rigid node] (2) [ right of=1047]      {$sl_2 \oplus \mathbb{C}$};

  \node (1050)[right of=105]{};  
\node (1051)[right of=1050]{};
\node (1052)[right of=1051]{};
\node (1053)[right of=1052]{};
\node (1054)[right of=1053]{};
\node (1054)[right of=1053]{};
\node (1055)[right of=1054]{};

 	\node[rigid node] (3) [ right of=1055]       {$g_5(\alpha)$};

\node (1060)[right of=106]{};  
\node (1061)[right of=1060]{};  
	\node[rigid node] (4) [ right of=1061]       {$g_4(\alpha,\beta)$};

	\node[rigid node] (5) [ right   of=107]      {$g_6$};
\node (1070)[right of=5]{};	
\node (1071)[right of=1070]{};
\node (10711)[right of=1071]{};
	\node[main node] (6) [ right of=10711]       {$n_4$};
\node (1072)[right of=6]{};
\node (1073)[right of=1072]{};
\node (1074)[right of=1073]{};
\node (1075)[right of=1074]{};

	\node[rigid node] (7) [right   of=1074]      {$g_3(\beta)$}; 

\node (1081)[right of=108]{};
\node (1080)[right of=1081]{};
	\node[main node] (8) [ right of=1080]      {$r_2 \oplus \mathbb{C}^2$};
\node (1082)[right of=8]{};
\node (1083)[right of=1082]{};
\node (1084)[right of=1083]{};
\node (1085)[right of=1084]{};
\node (1086)[right of=1085]{};
\node (1087)[right of=1086]{};
	\node[main node] (9) [ right      of=1086]       {$g_2(\beta)$};

\node (1100)[right of=110]{};
\node (1101)[right of=1100]{};	
\node (1102)[right of=1101]{};
\node (1103)[right of=1102]{};

	\node[main node] (10) [ right   of=1103]       {$n_3 \oplus \mathbb{C}$}; 
	
\node (1120)[right of=112]{};
\node (1121)[right of=1120]{};
\node (1122)[right of=1121]{};
\node (1123)[right of=1122]{};
\node (1124)[right of=1123]{};
\node (1125)[right of=1124]{};
\node (1126)[right of=1125]{};
\node (1127)[right of=1126]{};

	\node[main node] (11) [ right      of=1127]       {$g_1$};

\node (1160)[right of=116]{};
\node (1161)[right of=1160]{};
\node (1162)[right of=1161]{};
\node (1163)[right of=1162]{};
\node (1164)[right of=1163]{};
\node (1165)[right of=1164]{};
	\node[main node] (12) [ right  of=1165]       {$\mathbb{C}^4$};
 
\path[every node/.style={font=\sffamily\small}]

(1) edge node[above, fill=white] {\tiny $\alpha=0$} (3)
edge node[above, fill=white]{\tiny $\beta=0$} (4)

(2) edge node[above, fill=white]{\tiny $\alpha=-1$} (3)

(3) edge node[above, fill=white] {\tiny $\beta=\alpha+1$} (4)
 edge node[above, fill=white] {\tiny $\alpha=1,\beta=2$} (7)
 
 (4) edge node{} (6)
  edge node[above, fill=white] {\tiny $\alpha=\beta=0$} (8)
 edge node[above, fill=white] {\tiny $\alpha=1$} (9)
 
 (5) edge node{} (8)
 
 (6) edge node{} (10)
 
  (7) edge node{} (9)
  
 (8) edge node{} (10)
 
 (9) edge node{} (10)
 edge node[above, fill=white] {\tiny $\beta=1$} (11)
 
 (10) edge node{} (12)
 
 (11) edge node{} (12);
        
\end{tikzpicture}

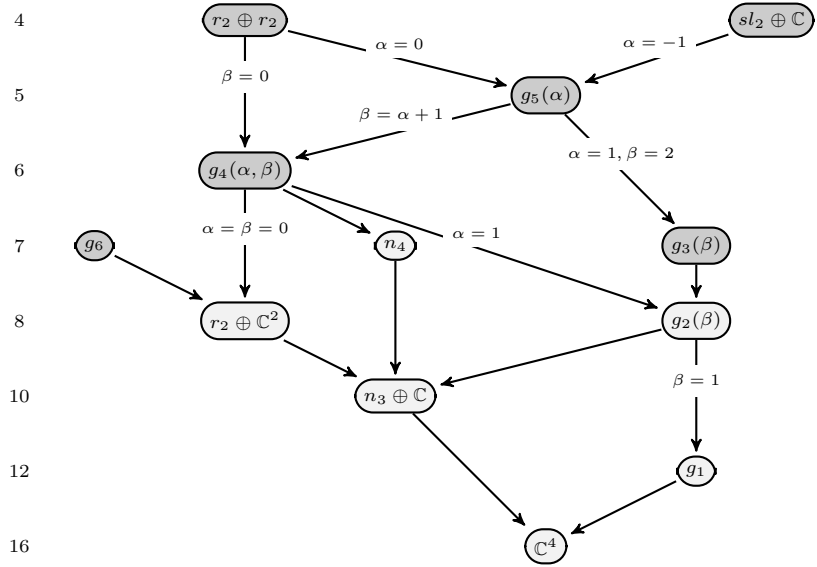
{\bf Figure I.}  The graph of primary degenerations for four-dimensional binary Lie algebras.

\end{center}

\begin{Proof} Tables 4.2 and 4.3 placed below give the proofs for all primary degenerations and non-degenerations including non-Lie algebras.
$$
\begin{array}{|c|c|}
\hline
\mbox{degenerations}  &  \mbox{parametrized bases}\\
\hline
\hline
g_3(\beta) \to g_2(\beta) & 
E_1^t=e_1+e_2, E_2^t=te_2, E_3^t=(1-\beta)e_3+e_4, E_4^t=e_3+e_4 \\
\hline
g_6 \to  r_2 \oplus \mathbb{C}^2 &
E_1^t=e_3, E_2^t=e_4, E_3^t=e_1, E_4^t=te_2 \\
\hline
\end{array}
$$
\begin{center} Table 4.2. Degenerations of binary Lie algebras of dimension 4.
\end{center}

$$
\begin{array}{|c|c|}
\hline
\mbox{non-degenerations}  &  \mbox{arguments}\\
\hline
\hline
g_3(\beta) \ \bcancel{\to} \  r_2 \oplus \mathbb{C}^2, g_1 (\beta\neq 1), g_2(\gamma \neq\beta )
&
\begin{array}{c}
c_{ij}(g_3(\beta)) = \frac{(\beta^i +2)(\beta^j+2)}{\beta^{i+j}+2}\mbox{, but }c_{ij}(r_2 \oplus \mathbb{C}^2)=1,
c_{ij}(g_1)=3\mbox{ and  }\\
c_{ij}(g_2(\gamma)) = \frac{(\gamma^i +2)(\gamma^j+2)}{\gamma^{i+j}+2}
\end{array}\\
\hline
g_6 \ \bcancel{\to} \ g_1, g_2(\beta) &  dim\,(g_6)^2< dim\,(g_2(\beta))^2 \leqslant dim\,(g_1)^2 \\ 
\hline
\end{array}
$$
\begin{center} Table 4.3. Non-degenerations of binary Lie algebras of dimension 4.
\end{center}
\end{Proof}

{\bf Remark.} Gorbatsevich classified all finite-dimensional anticommutative algebras of level 3 in \cite{gorb93}.
It follows  from Theorem \ref{first} that
his classification is not correct.
Namely, there are binary Lie algebras of level 3 in dimension 4, which are not included in the classification of Gorbatsevich.

\begin{corollary}\label{irr1}
$Irr(BL_4)=\{\mathcal{C}_i\}_{1\le i\le 5}$, where
$$
\begin{aligned}
\mathcal{C}_1&=\overline{O(sl_2 \oplus \mathbb{C})}=O\left(\{sl_2 \oplus \mathbb{C},g_5(-1),g_4(-1,0),n_4,n_3\oplus\mathbb{C},\mathbb{C}^4\}\right),\\
\mathcal{C}_2&=\overline{O(r_2\oplus r_2)}=O\left(\{r_2\oplus r_2,g_5(0),g_2(0), r_2 \oplus \mathbb{C}^2,n_4,n_3\oplus\mathbb{C},\mathbb{C}^4\}\cup\bigcup\limits_{\alpha\in\mathbb{C}}\{g_4(\alpha,0)\}\right),\\
\mathcal{C}_3&=\overline{\bigcup\limits_{\alpha\in\mathbb{C}} O(g_5(\alpha))}=O\left(\bigcup\limits_{\alpha\in\mathbb{C}}\{g_5(\alpha),g_4(\alpha,\alpha+1)\}\cup\{g_2(0),g_2(2),g_3(2),n_4,n_3\oplus\mathbb{C},\mathbb{C}^4\}\right),\\
\mathcal{C}_4&=\overline{\bigcup\limits_{\alpha,\beta\in\mathbb{C}} O(g_4(\alpha,\beta))}=O\left(\bigcup\limits_{\alpha,\beta\in\mathbb{C}}\{g_4(\alpha,\beta),g_2(\beta)\}\cup\{r_2\oplus\mathbb{C}^2,n_4,n_3\oplus\mathbb{C},g_1,\mathbb{C}^4\}\right),\\
\mathcal{C}_5&=\overline{\bigcup\limits_{\beta\in\mathbb{C}} O(g_3(\beta))}=O\left(\bigcup\limits_{\beta\in\mathbb{C}}\{g_3(\beta),g_2(\beta)\}\cup\{g_6,r_2\oplus\mathbb{C}^2,n_3\oplus\mathbb{C},g_1,\mathbb{C}^4\}\right).
\end{aligned}
$$
In particular, $Rig(BL_4)=Rig(Lie_4)=\{sl_2 \oplus \mathbb{C},r_2\oplus r_2\}$.
\end{corollary}
\begin{Proof} In view of Theorem \ref{first} and the fact that $Lie_4$ is a closed subset of $BL_4$ it is enough to prove that
\begin{multline*}
g_5(*)\not\to g_4(\alpha,\beta)\,(g_4(\alpha,\beta)\not\cong g_4(\gamma,\gamma+1)\mbox{ for any }\gamma\in\mathbb{C}),
g_5(*)\not\to r_2\oplus\mathbb{C}^2, g_5(*)\not\to g_2(\beta)\,(\beta\neq 0,2), g_5(*)\not\to g_1,\\
g_4(*,*)\not\to g_3(2),g_3(*)\to g_6, g_3(*)\not\to n_4,
\end{multline*}
where $g_5(*)=\{g_5(\alpha)\}_{\alpha\in\mathbb{C}}$, $g_4(*,*)=\{g_4(\alpha,\beta)\}_{\alpha,\beta\in\mathbb{C}}$ and $g_3(*)=\{g_3(\beta)\}_{\beta\in\mathbb{C}}$.
Let us define
$$
\mathcal{R}=\left\{A\left| \begin{array}{c}A=\langle f_1,f_2,f_3,f_4\rangle,\langle f_3,f_4\rangle^2=0,\langle f_2,f_3,f_4\rangle^2\subset\langle f_4\rangle,A\langle f_2,f_3,f_4\rangle\subset \langle f_2,f_3,f_4\rangle, A\langle f_3,f_4\rangle\subset\langle f_3,f_4\rangle,\\
A\langle f_4\rangle\subset\langle f_4\rangle,c_{1,2}^2+c_{1,3}^3=c_{1,4}^4,\mbox{ where } f_if_j=\sum\limits_{k=1}^4c_{i,j}^kf_k$ for all $1\leqslant i,j\leqslant 4\end{array}\right.\right\}.
$$
One can take $f_1=e_1$, $f_2=e_3$, $f_3=e_2$ and $f_4=e_4$ and check that $g_5(\alpha)\in\mathcal{R}$ for all $\alpha\in\mathbb{C}$. 

Let us prove that $g_2(\beta)\not\in\mathcal{R}$ if $\beta\not=0,2$. Assume that there is some basis $\tilde f_i$ ($1\le i\le 4$) of $V$ such that the structure constants $\tilde c_{i,j}^k$ of $g_2(\beta)$ in it satisfy all required conditions. Let $U=\langle \tilde f_2,\tilde f_3,\tilde f_4\rangle$ and $L:U\rightarrow U$ be the operator of left multiplication by $\tilde f_1$. It follows from the definition of $\mathcal{R}$ that the matrix of $L$ in the basis $\tilde f_2,\tilde f_3,\tilde f_4$ is lower triangular. Hence, $\tilde c_{1,2}^2$, $\tilde c_{1,3}^3$ and $\tilde c_{1,4}^4$ are eigen values of $L$. On the other hand, it is easy to see that $U =\langle e_2,e_3,e_4\rangle$ and $\tilde f_1= ce_1+v$ for some $c\in\mathbb{C}^*$ and $v\in U$. Then the eigen values of $L$ are $c$, $c$ and $\beta c$. Then we have $c=(\beta+1)c$ or $\beta c=2c$, i.e. $\beta=0$ or $\beta=2$.

Analogously one can prove that $g_4(\alpha,\beta)\not\in\mathcal{R}$ if $\alpha-\beta\not=1$, $\alpha-\beta\not=-1$ and $\alpha+\beta\not=1$, and $r_2\oplus\mathbb{C}^2,g_1\not\in\mathcal{R}$.

Since $\langle e_2,e_3,e_4\rangle$ is an abelian subalgebra of $g_4(\alpha,\beta)$ and there is no three-dimensional abelian subagebra in $g_3(2)$, we have $g_4(*,*)\not\to g_3(2)$ by Lemma \ref{gmain}. Let now define
$$
\mathcal{R}=\left\{A\left| \begin{array}{c}A=\langle f_1,f_2,f_3,f_4\rangle,\langle f_2,f_3,f_4\rangle^2\subset\langle f_4\rangle,A\langle f_2,f_3,f_4\rangle\subset \langle f_2,f_3,f_4\rangle, A\langle f_3,f_4\rangle\subset\langle f_3,f_4\rangle,
A\langle f_4\rangle\subset\langle f_4\rangle,\\c_{1,2}^2=c_{1,3}^3,c_{1,2}^3=0,\mbox{ where } f_if_j=\sum\limits_{k=1}^4c_{i,j}^kf_k$ for all $1\leqslant i,j\leqslant 4\end{array}\right.\right\}.
$$
One can take $f_i=e_i$ ($1\le i\le 4$) and check that $g_3(\beta)\in\mathcal{R}$ for all $\beta\in\mathbb{C}$. On the other hand, it is not hard to check that $n_4\not\in\mathcal{R}$. Finally, to prove that $g_3(*)\to g_6$ it is enough to take the parametrized basis
$$
E_1^t=e_2,E_2^t=e_3,E_3^t=e_4,E_4^t=-te_1
$$
and the parametrized index $\beta(t)=\frac{1}{t}$.
\end{Proof}

\begin{corollary}
$Irr(Mal_4)=\{\mathcal{C}_i\}_{1\le i\le 4}\cup\{\mathcal{C}_5'\}$, where $\mathcal{C}_i$ ($1\le i\le 4$) are the same as in Corollary \ref{irr1}, and $$\mathcal{C}_5'=\overline{O(g_3(-1))}=O\left(\{g_3(-1),g_2(-1),n_3\oplus\mathbb{C},\mathbb{C}^4\}\right).$$
In particular, $Rig(Mal_4)=\{sl_2 \oplus \mathbb{C},r_2\oplus r_2,g_3(-1)\}$.
\end{corollary}
\begin{Proof} Everything follows from Theorem \ref{first}, Corollary \ref{irr1} and Table 4.1.
\end{Proof}

\section{Degenerations of nilpotent Malcev  algebras of dimension 5}

For five-dimensional nilpotent Malcev algebras we have the following table, which is constructed using results of \cite{GRH} and \cite{K70}.

\begin{center}
\begin{equation*}
\begin{array}{|c|c|c|c|c|c|} 

\hline  A                       &\mbox{multiplication table}& Der(A)& Z_l(A) & A^2 & A^3  \\
\hline

\hline  n_3 \oplus \mathbb{C}^2   & \begin{array} {c} e_1e_2=e_3 \end{array} &  16   & 2+13 & 1 & 0     \\

\hline  n_4 \oplus \mathbb{C}   & \begin{array} {c} e_1e_2=e_3, e_1e_3=e_4 \end{array} &  12   & 1+124 & 2 & 1     \\

\hline  g_{5,1}   & \begin{array} {c} e_1e_2=e_5, e_3e_4=e_5 \end{array} &  11   & 15 & 1 & 0     \\

\hline  g_{5,2}   & \begin{array} {c} e_1e_2=e_4, e_1e_3=e_5  \end{array} &  15   & 25 & 2 & 0     \\

\hline  g_{5,3}   & \begin{array} {c} e_1e_2=e_3, e_1e_4=e_5,  e_2e_3=e_5 \end{array} &  10   & 135 & 2 & 1     \\

\hline  g_{5,4}   & \begin{array} {c} e_1e_2=e_3, e_1e_3=e_4, e_2e_3=e_5 \end{array} &  10   & 235 & 3 & 2     \\

\hline  g_{5,5}   & \begin{array} {c} e_1e_2=e_3, e_1e_3=e_4, e_1e_4=e_5 \end{array} &  9  & 1235 & 3 & 2    \\

\hline  g_{5,6}   & \begin{array} {c} e_1e_2=e_3, e_1e_3=e_4, e_1e_4=e_5, e_2e_3=e_5 \end{array} &  8   & 1235 & 3 & 2    \\

\hline
\hline  M_{5}   & \begin{array} {c} e_1e_2=e_4, e_3e_4=e_5 \end{array} &  9   & 135 & 2 & 0     \\

\hline
 \end{array}
\end{equation*}
\end{center}
\begin{center} Table 5.1. Nilpotent Malcev algebras of dimension 5.
\end{center}

\begin{Th}\label{graph5}
The graph of primary degenerations for nilpotent Malcev algebras of dimension 5 has the following form:
\end{Th}

\begin{center}

\begin{tikzpicture}[->,>=stealth',shorten >=0.05cm,auto,node distance=1.3cm,
                    thick,main node/.style={rectangle,draw,fill=gray!10,rounded corners=1.5ex,font=\sffamily \scriptsize \bfseries },rigid node/.style={rectangle,draw,fill=black!20,rounded corners=1.5ex,font=\sffamily \scriptsize \bfseries },style={draw,font=\sffamily \scriptsize \bfseries }]

  \node[rigid node] (1)                           {$g_{5,6}$ };

  \node[main node] (3) [ right       of=1]       {$g_{5,5}$};

  \node[main node] (2) [ right        of=3]       {$g_{5,4}$};

  \node[main node] (4) [ above        of=2]       {$g_{5,3}$};

  \node[rigid node] (55) [left          of=4]      {$M_5$};

  \node[main node] (8) [ right        of=2]       {$g_{5,1}$};

  \node (90) [ above         of=55]      {$9$};
  \node (91) [ left    of=90]      {$8$};
  \node (92) [ right    of=90]      {$10$};
  \node (93) [ right   of=92]      {$11$};
  \node (94) [ right   of=93]      {$12$};
  \node (94r) [right of=94]{};
  \node (95) [ right   of=94r]      {$15$};
  \node (95r) [right of=95]{};
  \node (96) [ right   of=95r]      {$16$};
  \node (96r) [right of=96]{};
  \node (97) [ right   of=96r]      {$25$};

  \node[main node] (5) [ below         of=94]       {$n_4 \oplus \mathbb{C}$};

\node (5r) [right of=5]{};

  \node[main node] (6) [ right       of=5r]       {$g_{5,2}$};
  
\node (6r) [right of=6]{};

  \node[main node] (7) [ right        of=6r]       {$n_3 \oplus \mathbb{C}^2$};
  
\node (7r) [right of=7]{};

  \node[main node] (111) [ right   of=7r]       {$\mathbb{C}^5$};
 
  \path[every node/.style={font=\sffamily\small}]
    (1) edge  [bend right] node[left] {} (2)
    (1) edge   node[left] {} (3)
    (1) edge   node[left] {} (4)

    (2) edge   node[left] {} (5)
    
    (3) edge   node[left] {} (5)
    
    (4) edge   node[left] {} (8)
    (4) edge     node[left] {} (5)
    
    (5) edge   node[left] {} (6)

    (6) edge   node[left] {} (7)

    (8) edge  node[left] {} (7)
    
    (7) edge   node[left] {} (111)

    (55) edge   node[left] {} (4);

\end{tikzpicture}

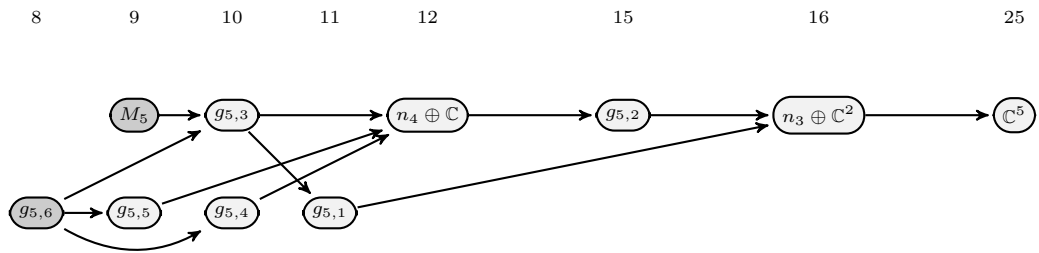
{\bf Figure II.} The graph of primary degenerations for five-dimensional nilpotent Malcev algebras.

\end{center}

\begin{Proof} It is enough to verify the assertions of the form $M_5\to A$ for such $A$ that $dim\,Der(A)<9$.
So we have to check that $M_5\to g_{5,3}$ and $M_5\not\to g_{5,4}$. The parametrized basis formed by $E_1^t=e_1-e_4$, $E_2^t=te_2+te_3$, $E_3^t=te_4+te_5$, $E_4^t=t^2e_3$ and $E_5^t=t^2e_5$ gives the required degeneration. 
The assertion $M_5\not\to g_{5,4}$ follows from the fact that $dim\,(g_{5,4})^2>dim\,(M_5)^2$.
\end{Proof}

\begin{corollary}
$Irr(NMal_5)=\{\mathcal{C}_1,\mathcal{C}_2\}$, where
$
\mathcal{C}_1=\overline{O(g_{5,6})}=NLie_5$ and $\mathcal{C}_2=\overline{O(M_5)}=NMal_5\setminus\{g_{5,6},g_{5,5},g_{5,4}\}.
$
In particular, $Rig(NMal_5)=\{g_{5,6},M_5\}$.
\end{corollary}
\begin{Proof}
Since there is only finite number of isomorphism classes of five-dimensional nilpotent Malcev algebras, everything  follows from Theorem \ref{graph5}.
\end{Proof}

\newpage
\section{Degenerations of nilpotent Malcev  algebras of dimension 6}

We use the table of invariants for nilpotent six-dimensional Lie algebras from \cite{S90} 
and classification of nilpotent six-dimensional Malcev non-Lie algebras from \cite{K70} to construct the table containing important invariants for nilpotent six-dimensional Malcev algebras. To simplify the notation we write $g_i$ instead of $g_{6,i}$, and $g_i^{\mathbb{C}}$ and $M_5^{\mathbb{C}}$ instead of $g_{5,i}\oplus \mathbb{C}$ and $M_5\oplus\mathbb{C}$ respectively.

\begin{center}
\begin{equation*}
\begin{array}{|c|c|c|c|c|c|} 

\hline  A                       &\mbox{multiplication table}&  Der(A)  & Z_l(A) & A^2  &  A^3 \\
\hline

\hline
\begin{array} {c}  \  g_{1}   \end{array}   
& \begin{array} {c} e_1e_2=e_3, e_1e_3=e_4, e_1e_4=e_6, 
e_2e_3=e_6, e_2e_5=e_6  \end{array} &   
\begin{array} {c} \ 11 \end{array}& 1346 & 3 & 2    \\

\hline
\begin{array} {c}  \  g_{2}   \end{array}   
& \begin{array} {c} e_1e_2=e_3, e_1e_3=e_4, 
e_1e_4=e_6, e_2e_5=e_6  \end{array} &   
\begin{array} {c} \ 12 \end{array}& 1346 & 3 & 2    \\

\hline
\begin{array} {c}  \  g_{3}   \end{array}   
& \begin{array} {c} e_1e_2=e_3, e_1e_3=e_6,  
e_4e_5=e_6 \end{array} &   
\begin{array} {c} \ 14 \end{array}& 146 & 2 & 1    \\

\hline
\begin{array} {c}  \  g_{5}   \end{array}   
& \begin{array} {c} e_1e_2=e_3, e_1e_3=e_4, e_1e_4=e_5, 
e_1e_5=e_6, e_2e_3=e_5, e_2e_4=e_6 \end{array} &   
\begin{array} {c} \ 9 \end{array}& 12346 & 4 & 3    \\

\hline
\begin{array} {c}  \  g_{6}   \end{array}   
& \begin{array} {c} e_1e_2=e_3, e_1e_3=e_4, e_1e_4=e_5, 
e_2e_3=e_5, e_2e_5=e_6, e_3e_4=-e_6 \end{array} &   
\begin{array} {c} \ 8 \end{array}& 12346 & 4 & 3    \\

\hline
\begin{array} {c}  \  g_{7}   \end{array}   
& \begin{array} {c} e_1e_2=e_3, e_1e_3=e_4, e_1e_4=e_5, 
e_1e_5=e_6, e_2e_3=e_6  \end{array} &   
\begin{array} {c} \ 10 \end{array}& 12346 & 4 & 3    \\

\hline
\begin{array} {c}  \  g_{8}   \end{array}   
& \begin{array} {c} e_1e_2=e_3, e_1e_3=e_4, e_2e_5=e_6, 
e_3e_4=-e_6  \end{array} &   
\begin{array} {c} \ 9 \end{array}& 12346 & 3 & 2    \\

\hline
\begin{array} {c}  \  g_{9}   \end{array}   
& \begin{array} {c} e_1e_2=e_3, e_1e_3=e_4, e_1e_4=e_5, 
e_1e_5=e_6  \end{array} &   
\begin{array} {c} \ 11 \end{array}& 12346 & 4 & 3    \\

\hline
\begin{array} {c}  \  g_{10}   \end{array}   
& \begin{array} {c} e_1e_2=e_4, e_1e_3=e_5,  
e_1e_4=e_6, e_3e_5=e_6  \end{array} &  
\begin{array} {c} \ 12 \end{array}& 136 & 3 & 1    \\

\hline
\begin{array} {c}  \  g_{14}   \end{array}   
& \begin{array} {c} e_1e_2=e_3, e_1e_3=e_4, e_1e_5=e_6, 
e_2e_3=e_5, e_2e_4=e_6  \end{array} &   
\begin{array} {c} \ 10 \end{array}& 1346 & 4 & 3    \\

\hline
\begin{array} {c}  \  g_{15}   \end{array}   
& \begin{array} {c} e_1e_2=e_3, e_1e_3=e_5,   
e_1e_4=e_6, e_2e_3=e_6  \end{array} &   
\begin{array} {c} \ 13 \end{array}& 246 & 3 & 2    \\


\hline
\begin{array} {c}  \  g_{16}   \end{array}   
& \begin{array} {c} e_1e_2=e_3, e_1e_3=e_5,   
e_1e_4=e_5, e_2e_3=e_6  \end{array} &   
\begin{array} {c} \ 12 \end{array}& 246 & 3 & 2    \\

\hline
\begin{array} {c}  \  g_{17}   \end{array}   
& \begin{array} {c} e_1e_2=e_3, e_1e_3=e_5,   
e_1e_4=e_6  \end{array} &   
\begin{array} {c} \ 15 \end{array}& 246 & 3 & 1    \\

\hline
\begin{array} {c}  \  g_{18}   \end{array}   
& \begin{array} {c} e_1e_2=e_3, e_1e_3=e_5,   
e_2e_4=e_6  \end{array} &   
\begin{array} {c} \ 13 \end{array}& 246 & 3 & 1    \\

\hline
\begin{array} {c}  \  g_{20}   \end{array}   
& \begin{array} {c} e_1e_2=e_3, e_1e_3=e_5,  
e_1e_4=e_6, e_2e_4=e_5  \end{array} &   
\begin{array} {c} \ 14 \end{array}& 246 & 3 & 1    \\

\hline
\begin{array} {c}  \  g_{21}   \end{array}   
& \begin{array} {c} e_1e_2=e_5, e_1e_3=e_6,  
e_3e_4=e_5  \end{array} &  
\begin{array} {c} \ 17 \end{array}& 26 & 2 & 0    \\

\hline
\begin{array} {c}  \  g_{23}   \end{array}   
& \begin{array} {c} e_1e_2=e_3, e_1e_3=e_4,   
e_1e_4=e_5, e_2e_3=e_6  \end{array} &   
\begin{array} {c} \ 11 \end{array}& 2346 & 4 & 3    \\

\hline
\begin{array} {c}  \  g_{24}   \end{array}   
& \begin{array} {c} e_1e_2=e_4, e_1e_3=e_5,   
e_2e_3=e_6  \end{array} &   
\begin{array} {c} \ 16 \end{array}& 36 & 3 & 0    \\

\hline
\begin{array} {c}  \  g_{1}^\mathbb{C}   \end{array}   
& \begin{array} {c} e_1e_2=e_5, e_3e_4=e_5  
  \end{array} &   
\begin{array} {c} \ 21 \end{array}& 1+15 & 1 & 0    \\

\hline

\begin{array} {c}  \  g_{2}^\mathbb{C}   \end{array}   
& \begin{array} {c} e_1e_2=e_4, e_1e_3=e_5   \end{array} &   
\begin{array} {c} \ 19 \end{array}& 1+25 & 2 & 0    \\

\hline
\begin{array} {c}  \  g_{3}^\mathbb{C}   \end{array}   
& \begin{array} {c} e_1e_2=e_3, e_1e_4=e_5,  
e_2e_3=e_5  \end{array} &   
\begin{array} {c} \ 15 \end{array}& 1+135 & 2 & 1    \\

\hline
\begin{array} {c}  \  g_{4}^\mathbb{C}   \end{array}   
& \begin{array} {c} e_1e_2=e_3, e_1e_3=e_4,  
e_2e_3=e_5  \end{array} &   
\begin{array} {c} \ 15 \end{array}& 1+235 & 3 & 2    \\

\hline
\begin{array} {c}  \  g_{5}^\mathbb{C}   \end{array}   
& \begin{array} {c} e_1e_2=e_3, e_1e_3=e_4,   
e_1e_4=e_5    \end{array} &   
\begin{array} {c} \ 13 \end{array}& 1+1235 & 3 & 2    \\

\hline
\begin{array} {c}  \  g_{6}^\mathbb{C}   \end{array}   
& \begin{array} {c} e_1e_2=e_3, e_1e_3=e_4,  
e_1e_4=e_5, e_2e_3=e_5  \end{array} &   
\begin{array} {c} \ 12
\end{array}& 1+1235 & 3 & 2   \\

\hline

\begin{array} {c}  \  n_3 \oplus n_3  \end{array}   
& \begin{array} {c} e_1e_3=e_5, e_2e_4=e_6  
   \end{array} &   
\begin{array} {c} \ 16 \end{array}& 13+13 & 2 & 0    \\

\hline
\begin{array} {c}  \  n_4 \oplus \mathbb{C}^2   \end{array}   
& \begin{array} {c} e_1e_2=e_3, e_1e_3=e_4  \end{array} &   
\begin{array} {c} \ 17 \end{array}& 2+124 & 2 & 1    \\

\hline
\begin{array} {c}  \  n_3 \oplus \mathbb{C}^3  \end{array}   
& \begin{array} {c} e_1e_2=e_3 \end{array} &   
\begin{array} {c} \ 24 \end{array}& 3+13 & 1 & 0    \\

\hline
\hline
\begin{array} {c}  \  M_5^\mathbb{C}  \end{array}   
& \begin{array} {c} e_1e_2=e_5, e_3e_5=e_6  \end{array} &   
\begin{array} {c} \ 14 \end{array}& 1+135 & 2 & 1    \\

\hline \begin{array} {c}  \ M_{1}^{0, 1} \end{array}   
& \begin{array} {c} e_1e_2=e_5, e_3e_4=  e_5, e_3e_5=e_6 \end{array} &   
\begin{array} {c} \ 13  \end{array}& 126 &2 & 1    \\

\hline \begin{array} {c}  \ M_{1}^{1, 0} \end{array}   
& \begin{array} {c} e_1e_2=e_5,  e_1e_4=  e_6, e_3e_5=e_6 \end{array} &   
\begin{array} {c} \ 12  \end{array}& 136  &2  & 1   \\

\hline \begin{array} {c}  \ M_{1}^{1, 1}  \end{array}   
& \begin{array} {c} e_1e_2=e_5,   e_1e_4=  e_6, e_3e_4=  e_5, e_3e_5=e_6 \end{array} &   
\begin{array} {c} \ 11 \end{array}& 126&2 & 1    \\

\hline  \begin{array} {c} \ M_{2}^0  \end{array}   
& \begin{array} {c} e_1e_2=e_4, e_1e_3=e_5, e_2e_5=  e_6  \end{array} & 12     & 246 &3& 1  \\

\hline   M_{2}^{-1} 
& \begin{array} {c} e_1e_2=e_4, e_1e_3=e_5,  
e_2e_5= e_6, e_3e_4= -e_6 \end{array} &  13 & 136 & 3 & 1  \\

\hline  \begin{array} {c} \ M_{2}^{\epsilon}, \epsilon \neq -1,0 \end{array}   
& \begin{array} {c} e_1e_2=e_4, e_1e_3=e_5,   
e_2e_5= e_6, e_3e_4= \epsilon e_6 \end{array} &  11 & 136 & 3 & 1  \\

\hline  \begin{array} {c} \ M_{3} \end{array}   
& \begin{array} {c} e_1e_2=e_4, e_1e_3=e_5,   
e_2e_4=e_6, e_2e_5= e_6, e_3e_4=-e_6 \end{array} &  11 & 136 & 3 & 1  \\

\hline  \begin{array} {c} \ M_{4} \end{array}    
& \begin{array} {c} e_1e_2=e_4, e_1e_3=e_5,   e_1e_5=e_6,  e_3e_4= e_6,  \end{array} &  13  & 136 & 3 & 1    \\

\hline \ M_{5}^{0}
 & \begin{array} {c} e_1e_2=e_4, e_2e_4= e_5,   e_3e_4=  e_6 \end{array} &  
  11 &246  & 3 &  2 \\

\hline \ M_{5}^{1}
 & \begin{array} {c} e_1e_2=e_4, e_1e_3=  e_5, e_2e_4= e_5,   e_3e_4=  e_6 \end{array} &  
   10  &246  & 3 &  2 \\

\hline  \ M_{6}^{0}    
& \begin{array} {c} e_1e_2=e_3,   e_1e_3=  e_5, e_1e_5=e_6, e_3e_4=e_6 \end{array} & 
 10 &  1346  & 3 &  2  \\

\hline  \ M_{6}^{\epsilon}, \epsilon \neq 0     
& \begin{array} {c} e_1e_2=e_3,    e_1e_3=  e_5, e_1e_5=e_6, e_2e_4=\epsilon e_5, e_3e_4=e_6 \end{array} & 
 10 & 1236  & 3 &  2  \\

\hline  \ M_{7}^{0}
& \begin{array} {c} e_1e_2=e_4, e_1e_4 =  e_5,  e_1e_5=  e_6, e_2e_3=e_5  \end{array} & 
   11 
& 1246 & 3 & 2    \\

\hline   \ M_{7}^{1} 
& \begin{array} {c} e_1e_2=e_4,  e_1e_4 =  e_5,   e_1e_5=  e_6, e_2e_3=e_5, e_2e_4=  e_6  \end{array} & 
 \begin{array} {c}   10   \end{array}   
& 1246 & 3 & 2     \\

\hline
 \end{array}
\end{equation*}
\end{center}
\begin{center} Table 6.1. Nilpotent Malcev algebras of dimension 6.
\end{center}

The algebra $M_2^{\epsilon}$ is isomorphic to $M_2^{\epsilon'}$ iff
$\epsilon\epsilon'=1$ or $\epsilon=\epsilon'.$ Apart from this exception any two algebras with different names from Table 6.1 are not isomorphic.

\newpage
\begin{Th}\label{third}
The graph of primary degenerations for nilpotent Malcev algebras of dimension 6 has the form presented in Figure III.
\end{Th}

\begin{Proof} Tables 6.2 and 6.3 presented below give the proofs for all primary degenerations and non-degenerations including non-Lie algebras.

\begin{center}
$$\begin{array}{|c|c|}
\hline
\mbox{degenerations}  &  \mbox{parametrized bases}\\
\hline
\hline
M_7^1 \to  g_{1}  &
\begin{array}{c}E_1^t=e_1-\frac{e_2}{t}-\frac{e_3}{t},E_2^t=e_2,E_3^t=e_4+\frac{e_5}{t},E_4^t=e_5,\\E_5^t=te_3+e_4+\frac{e_5}{t},E_6^t=e_6\end{array}\\\hline
M_7^1 \to  M_7^0 & E_1^t=e_1,E_2^t=te_2,E_3^t=e_3,E_4^t=te_4,E_5^t=te_5,E_6^t=te_6 \\\hline
M_7^1 \to M_4 &E_1^t=e_2, E_2^t=-e_3, E_3^t=-te_1, E_4^t=-e_5, E_5^t=te_4,E_6^t=te_6 \\\hline
M_7^1 \to M^{1,1}_1 & E_1^t=e_2, E_2^t=\frac{e_3}{t}, E_3^t=e_1, E_4^t=\frac{e_4}{t}, E_5^t=\frac{e_5}{t},E_6^t=\frac{e_6}{t}
\\
\hline
\hline
M_5^1 \to M_5^0 & E_1^t=e_1,E_2^t=e_2,E_3^t=te_3,E_4^t=e_4,E_5^t=e_5,E_6^t=te_6\\
\hline
\hline 
M_6^{\epsilon} \to g_1 &\begin{array}{c}
E_1^t=e_1+\frac{e_4}{2t},E_2^t=e_2+\frac{(\epsilon+2)e_3}{2t}-e_4,E_3^t=e_3+\frac{e_5}{t}-\frac{(\epsilon+2)e_6}{4t^2},\\E_4^t=e_5+\frac{e_6}{2t},E_5^t=te_4,E_6^t=e_6\end{array}\\\hline
M_6^{\epsilon} \to M_2^{\epsilon}\,\,(\epsilon\neq0) &
E_1^t=e_2,E_2^t=\frac{e_4}{t},E_3^t=e_1,E_4^t=\frac{\epsilon e_5}{t},E_5^t=-e_3,E_6^t=\frac{e_6}{t}\\\hline
M_6^{\epsilon} \to M_5^0\,\,(\epsilon=0)&
E_1^t=-e_2,E_2^t=e_1,E_3^t=-\frac{e_4}{t},E_4^t=e_3,E_5^t=e_5,E_6^t=\frac{e_6}{t} \\\hline
M_6^{\epsilon} \to M_1^{1,0}\,\,(\epsilon=0) &
E_1^t=te_1,E_2^t=e_2,E_3^t=-e_4,E_4^t=e_5,E_5^t=te_3,E_6^t=te_6 \\\hline
M_6^{\epsilon} \to  M_4 \,\,(\epsilon=0) &
E_1^t=e_1,E_2^t=e_2,E_3^t=\frac{e_3}{t}-\frac{e_4}{t},E_4^t=e_4,E_5^t=\frac{e_5}{t},E_6^t=\frac{e_6}{t} \\\hline
M_6^{\epsilon} \to M_3\,\,(\epsilon=-1) &
\begin{array}{c}E_1^t=-e_2-\frac{e_3}{2},E_2^t=te_1,E_3^t=te_4,E_4^t=te_3+\frac{te_5}{2},\\E_5^t=te_5-\frac{te_6}{2},E_6^t=t^2e_6\end{array}\\
\hline
\hline
M_2^{\epsilon} \to g_{10}\,\,(\epsilon\neq-1,0) &\begin{array}{c}E_1^t=\frac{\epsilon e_1}{t}+e_2-e_3,E_2^t=te_1,E_3^t=\frac{(\epsilon+1)te_3}{\epsilon}-e_1,E_4^t=(\epsilon+1)te_5,\\E_5^t=e_4+\epsilon e_5,E_6^t=(\epsilon+1)te_6 \\\end{array}\\ 
\hline
M_2^{\epsilon} \to  g_{18}\,\,(\epsilon=0) & E_1^t=e_2+e_3,E_2^t=e_1,E_3^t=-e_4-e_5,E_4^t=te_3,E_5^t=-e_6,E_6^t=te_5 \\ 
\hline
M_2^{\epsilon} \to M_5^{\mathbb{C}}\,\,(\epsilon=0) & 
E_1^t=e_1,E_2^t=e_3,E_3^t=e_2,E_4^t=\frac{e_4}{t},E_5^t=e_5,E_6^t=e_6 \\ 
\hline
M_2^{\epsilon} \to  g_{20}\,\,(\epsilon=-1) &  E_1^t=e_1-e_2,E_2^t=e_3,E_3^t=e_5,E_4^t=te_2+e_4,E_5^t=-e_6,E_6^t=te_4 \\
\hline
\hline
M_7^0 \to  g_5^{\mathbb{C}} & E_1^t=e_1,E_2^t=e_2,E_3^t=e_4,E_4^t=e_5,E_5^t=e_6,E_6^t=te_3 \\
\hline
M_7^0 \to g_3  & E_1^t=te_1,E_2^t=\frac{e_4}{t^2}-\frac{e_5}{t^2},E_3^t=\frac{e_5}{t}-\frac{e_6}{t},E_4^t=e_2,E_5^t=e_3,E_6^t=e_6  \\
\hline
M_7^0 \to M_2^0 & E_1^t=e_2, E_2^t=e_1, E_3^t=\frac{e_3}{t}, E_4^t=-e_4, E_5^t=\frac{e_5}{t}, E_6^t=\frac{e_6}{t} \\
\hline 
M_7^0 \to M_1^{0,1} &  E_1^t=te_3, E_2^t=-\frac{e_2}{t}, E_3^t=te_1, E_4^t=\frac{e_4}{t}, E_5^t=e_5, E_6^t=t e_6 \\
\hline
\hline 
M_5^0 \to g_{16} &  E_1^t=te_2,E_2^t=-e_1-e_3-e_4,E_3^t=te_4-te_5,E_4^t=te_3+te_4, E_5^t=t^2e_5,E_6^t=-te_6 \\
\hline 
M_5^0 \to M_2^0  & E_1^t=e_1+e_2,E_2^t=e_3+e_4,E_3^t=te_2,E_4^t=e_5,E_5^t=te_4,E_6^t=te_6\\
\hline
\hline 
M_1^{1,1} \to M_1^{0,1} & E_1^t=te_1,E_2^t=e_2,E_3^t=e_3,E_4^t=te_4,E_5^t=te_5,E_6^t=te_6 \\
\hline 
M_1^{1,1} \to  M_1^{1,0} & E_1^t=e_1,E_2^t=e_2,E_3^t=te_3,E_4^t=te_4,E_5^t=e_5,E_6^t=te_6
\\
\hline
\hline
M_3 \to g_{10} & E_1^t=e_1+e_2+e_3,E_2^t=t^2e_2+t^2e_3,E_3^t=te_2,E_4^t=t^2e_4+t^2e_5, E_5^t=te_4,E_6^t=t^2e_6  \\
\hline
M_3 \to M_2^{-1} & E_1^t=e_1,E_2^t=te_2,E_3^t=e_3,E_4^t=te_4,E_5^t=e_5,E_6^t=te_6 \\
\hline
\hline
M_1^{1,0} \to g_3 & E_1^t=e_1+e_3,E_2^t=te_2,E_3^t=te_5,E_4^t=te_3,E_5^t=-e_4+e_5,E_6^t=te_6 \\
\hline
\hline
M_1^{1,0} \to M_5^{\mathbb{C}} & E_1^t=e_1,E_2^t=e_2,E_3^t=e_3,E_4^t=te_4,E_5^t=e_5,E_6^t=e_6 \\
\hline
\hline
M_1^{0,1} \to g_3 & E_1^t=te_3,E_2^t=\frac{e_4}{t^2}-\frac{e_5}{t^2},E_3^t=\frac{e_5}{t}-\frac{e_6}{t},E_4^t=e_1,E_5^t=e_2,E_6^t=e_6 \\
\hline
M_1^{0,1} \to M_5^{\mathbb{C}}& E_1^t=e_1,E_2^t=e_2,E_3^t=e_3,E_4^t=te_4,E_5^t=e_5,E_6^t=e_6\\
\hline
\hline
M_4 \to g_{20} &  E_1^t=\frac{e_1}{t},E_2^t=e_3,E_3^t=\frac{e_5}{t},E_4^t=e_2+\frac{e_4}{t^2},E_5^t=\frac{e_6}{t^2},E_6^t=\frac{e_4}{t} \\
\hline
M_4 \to g_3 & E_1^t=e_1+e_3,E_2^t=e_2,E_3^t=e_4,E_4^t=te_3,E_5^t=\frac{e_4}{t}-\frac{e_5}{t},E_6^t=e_6 \\
\hline
M_4 \to  M_5^{\mathbb{C}} & E_1^t=te_1,E_2^t=\frac{e_2}{t},E_3^t=e_3,E_4^t=e_5,E_5^t=e_4,E_6^t=e_6 \\
\hline
\hline
M_5^{\mathbb{C}} \to g_3^{\mathbb{C}} & E_1^t=e_1-e_4, E_2^t=te_2+te_3, E_3^t=te_4+te_5, E_4^t=t^2e_3, E_5^t=t^2e_5,E_6^t=e_6 \\
\hline 
M_5^{\mathbb{C}} \to n_3\oplus n_3 & E_1^t=e_1,E_2^t=e_3,E_3^t=te_2,E_4^t=e_4+e_5,E_5^t=te_5,E_6^t=e_6 \\
\hline
\end{array}$$
\end{center}
\begin{center} Table 6.2. Degenerations of nilpotent Malcev algebras of dimension 6.
\end{center}

$$
\begin{array}{|c|c|}
\hline
 \mbox{non-degenerations}  &  \mbox{arguments}\\
\hline \hline
M_7^1 \ \bcancel{\to} \ g_9,g_{23} &  dim\,(g_9)^2=dim\,(g_{23})^2>dim\,(M_7^1)^2 \\
\hline
M_7^1 \  \bcancel{\to} \ M_2^{\epsilon}\,\,(\epsilon\neq0),M_5^0  &
\begin{array}{c}
\mathcal{R}=\left\{A\left| \begin{array}{c}A=\langle f_1,f_2,f_3,f_4,f_5,f_6\rangle,A^2\subset \langle f_4,f_5,f_6\rangle,\\
						\langle f_3,f_4,f_5,f_6\rangle^2=0,A\langle f_3,f_4,f_5,f_6\rangle\subset\langle f_5,f_6\rangle\end{array}\right.\right\}\vspace{0.1cm}\\
M_7^1\in\mathcal{R} \mbox{ (take $f_i=e_i$ for $1\leqslant i\leqslant 6$), but }M_2^{\epsilon},M_5^0\not\in\mathcal{R}
\end{array}\\
\hline
\hline
M_5^1 \ \bcancel{\to} \ \left\{\begin{array}{c}M_2^{\epsilon}\,\,(\epsilon\neq0),M_7^0,M_1^{1,0},\\M_1^{0,1},M_4,g_3 \end{array} \right\}&
\begin{array}{c}dim\,Z(M_2^{\epsilon})=dim\,Z(M_7^0)=dim\,Z(M_1^{1,0})\\=dim\,Z(M_1^{0,1})=dim\,Z(M_4)=dim\,Z(g_3)<dim\,Z(M_5^1) \end{array}\\
\hline
M_5^1 \ \bcancel{\to} \ g_5^{\mathbb{C}} &
dim\,Z_2(g_5^{\mathbb{C}})<dim\,Z_2(M_5^1)\\
\hline
\hline
M_6^{\epsilon} \ \bcancel{\to} \  g_9,g_{23} & dim\,(g_9)^2=dim\,(g_{23})^2>dim\,(M_6^{\epsilon})^2 \\
\hline
M_6^{\epsilon} \ \bcancel{\to} \  M_1^{0,1}\,\,(\epsilon=0) & dim\,Z_2(M_1^{0,1})<dim\,Z_2(M_6^0) \\
\hline
M_6^{\epsilon} \ \bcancel{\to} \ \left\{ \begin{array}{c} M_2^{\epsilon'}\,\,(M_2^{\epsilon'}\not\cong M_2^{\epsilon}),\\M_5^{\mathbb{C}}\,\,(\epsilon\neq0)   \end{array} \right\}&
\begin{array}{c} 
\mathcal{R}=\left\{A\left| \begin{array}{c}A=\langle f_1,f_2,f_3,f_4,f_5,f_6\rangle,A^2\subset \langle f_4,f_5,f_6\rangle,\\
						\langle f_3,f_4,f_5,f_6\rangle^2=0,\\x(yz)=\epsilon y(xz)\forall x\in A, y,z\in\langle f_2,f_3,f_4,f_5,f_6\rangle\end{array}\right.\right\}\vspace{0.1cm}\\
M_6^{\epsilon}\in\mathcal{R} \mbox{ (take $f_1=e_1$, $f_2=e_4$, $f_3=e_2$, $f_4=e_3$,$f_5=e_5$}\\
\mbox{ and $f_6=e_6$), but }M_2^{\epsilon'}\not\in\mathcal{R}\mbox{ and, if $\epsilon\neq 0$, then }M_5^{\mathbb{C}}\not\in\mathcal{R}
\end{array}\\
\hline
\hline
M_2^{\epsilon} \ \bcancel{\to} \ g_{5}^{\mathbb{C}},g_{4}^{\mathbb{C}}& dim\,(g_{5}^{\mathbb{C}})^3=dim\,(g_{4}^{\mathbb{C}})^3>dim\,(M_2^{\epsilon})^3 \\
\hline
M_2^{\epsilon} \ \bcancel{\to} \ n_3\oplus n_3 (\epsilon=-1) & 
\begin{array}{c}
\mathcal{R}=\left\{A\left| \begin{array}{c}A=\langle f_1,f_2,f_3,f_4,f_5,f_6\rangle,A^2\subset \langle f_4,f_5,f_6\rangle,\\
						\langle f_3,f_4,f_5,f_6\rangle^2=0,A\langle f_4,f_5,f_6\rangle\subset\langle f_6\rangle,\\
						\langle f_2, f_3,f_4,f_5,f_6\rangle\langle f_5,f_6\rangle=0,\\
						\langle f_2, f_3,f_4,f_5,f_6\rangle\langle f_3,f_4,f_5,f_6\rangle\subset\langle f_5,f_6\rangle,\\
						c_{2,3}^{5}c_{1,5}^6+c_{2,4}^6c_{1,3}^4=0,c_{2,3}^{5}c_{1,4}^6=c_{2,4}^6c_{1,3}^5, \\ \mbox{ where } f_if_j=\sum\limits_{k=1}^6c_{i,j}^kf_k$ for all $1\leqslant i,j\leqslant 6\end{array}\right.\right\}\vspace{0.1cm}\\
M_2^{-1}\in\mathcal{R} \mbox{ (take $f_1=e_2$, $f_2=e_3$, $f_3=e_1$, $f_4=e_4$,$f_5=e_5$}\\
\mbox{ and $f_6=e_6$), but }n_3\oplus n_3\not\in\mathcal{R}
\end{array}\\
\hline
\hline
M_7^0 \ \bcancel{\to} \ g_{10},M_1^{1,0},M_4,g_4^{\mathbb{C}} &
\begin{array}{c}
\mathcal{R}=\left\{A\left| \begin{array}{c}A=\langle f_1,f_2,f_3,f_4,f_5,f_6\rangle,A^2\subset \langle f_4,f_5,f_6\rangle,\\
						\langle f_2,f_3,f_4,f_5,f_6\rangle\langle f_4,f_5,f_6\rangle=0\end{array}\right.\right\}\vspace{0.1cm}\\
M_7^0\in\mathcal{R} \mbox{ (take $f_i=e_i$ for $1\leqslant i\leqslant 6$), but }g_{10},M_1^{1,0},M_4,g_4^{\mathbb{C}}\not\in\mathcal{R}
\end{array}\\
\hline
\hline
M_1^{1,1} \ \bcancel{\to} \  M_4,g_{17},g_{24} &
\begin{array}{c} dim\,(M_4)^2=dim\,(g_{17})^2=dim\,(g_{24})^2>dim\,(M_1^{1,1})^2
\end{array}\\
\hline
\hline
M_3 \ \bcancel{\to} \  g_{5}^{\mathbb{C}},g_{4}^{\mathbb{C}} &
dim\,(g_{5}^{\mathbb{C}})^3=dim\,(g_{4}^{\mathbb{C}})^3>dim\,(M_3)^3\\
\hline
\hline
M_1^{1,0} \ \bcancel{\to} \ M_1^{0,1} & dim\,Z_2(M_1^{0,1})<dim\,Z_2(M_1^{1,0})\\
\hline
\hline
M_4 \ \bcancel{\to} \  g_4^{\mathbb{C}} &
dim\,(g_{4}^{\mathbb{C}})^3>dim\,(M_4)^3 \\
\hline
\end{array}
$$
\begin{center} Table 6.3. Non-degenerations of nilpotent Malcev algebras of dimension 6.
\end{center}
\end{Proof}

\begin{center}
\begin{tikzpicture}[->,>=stealth',shorten >=0.05cm,auto,node distance=1.3cm,
                    thick,main node/.style={rectangle,draw,fill=gray!10,rounded corners=1.5ex,font=\sffamily \tiny \bfseries },rigid node/.style={rectangle,draw,fill=black!20,rounded corners=1.5ex,font=\sffamily \tiny \bfseries },style={draw,font=\sffamily \scriptsize \bfseries }]

  \node[main node] (11)        { $g_{6}^{\mathbb{C}}$};
  
    \node[main node] (12) [ right      of=11]       { $g_{16}$};

  \node[main node] (14) [ right      of=12]       { $g_{2}$};

  \node[main node] (13) [ right       of=14]       { $g_{10}$};
  
  \node (13r) [right of=13] {};
  
    \node[main node] (m20) [ right       of=13r]       { $M_2^0$};
    
    \node (m20r) [right of=m20] {};
    
    \node (m20rr) [right of=m20r] {};
    
    \node (m20rrr) [right of=m20rr] {};
    
      \node[main node] (m110) [ right       of=m20rrr]       { $M_1^{1, 0}$};

  \node (11a) [above of=11] {};
  
    \node[main node] (7) [ above      of=11a]       { $g_{9}$};

  \node[main node] (8) [ right       of=7]       { $g_{23}$};

  \node[main node] (10) [ right       of=8]       { $g_{1}$};
  
    \node[main node] (9) [ right       of=10]       { $M_2^{\epsilon}$};
    
    \node[main node] (m70) [  right      of=9]       { $M_7^0$};
    
    \node[main node] (m50) [  right      of=m70]       { $M_5^0$};
    
    \node (m50r) [right of=m50] {};
    
    \node[main node] (m111) [  right       of=m50r]       { $M_1^{1,1}$};
    
    \node (m111r) [right of=m111] {};
    
    \node[main node] (m3) [right of=m111r] {$M_3$};

  \node (7a) [above of=7] {};
  
  \node[main node] (4) [above       of=7a]       { $g_{7}$};

  \node[main node] (5) [  right      of=4]       { $g_{14}$};
  
    \node (5r) [  right      of=5]       {};
    
    \node (5rr) [  right      of=5r]       {};
    
    \node[rigid node] (m71) [  right      of=5rr]       { $M_7^1$};
    
    \node[rigid node] (m51) [  right      of=m71]       { $M_5^1$};
    
    \node (m51r) [right of=m51] {};
    
    \node (m51rr) [right of=m51r] {};

  \node[rigid node] (6) [  right      of=m51rr]       { $M_6^{\epsilon}$};

\node (4a) [above left of=4] {};

  \node(4a1) [above right       of=4a]       {};  

  \node[main node] (3) [right       of=4a1]       { $g_{5}$};  
  
  \node (3r) [right of=3] {};

        \node[main node] (2) [ right      of=3r]       { $g_{8}$};

  \node[rigid node] (1) [above  of=3r]                          { $g_{6}$ };

\node (12b) [below left of=12] {};
 
  \node[main node] (15) [ below right       of=12b]       { $g_{5}^{\mathbb{C}}$};

  \node[main node] (16) [  right       of=15]       { $g_{15}$};
  
  \node (16r) [right of=16] {};

  \node[main node] (17) [ right       of=16r]       { $g_{18}$};
  
      \node[main node] (m2-1) [ right       of=17]       { $M_2^{-1}$};
      
      \node (m2-1r) [right of=m2-1] {};
      
      \node[main node] (m101) [  right       of=m2-1r]       { $M_1^{0,1}$};
      
      \node (m101r) [right of=m101] {};
      
      \node[main node] (m4) [right of=m101r] {$M_4$};

\node (16b) [below left of=16] {}; 

  \node[main node] (18) [ below right       of=16b]       { $g_{20}$};
  
  \node (18r) [right of=18] {};

  \node[main node] (19) [  right       of=18r]       { $g_{3}$};
  
    \node (19r) [ right       of=19]       {};
  
    \node (19rr) [ right       of=19r]       {};

  \node[main node] (m5c) [ right       of=19rr]       { $M_5^{\mathbb{C}}$};

\node (18b) [below left of=18] {}; 

  \node[main node] (20) [ below right        of=18b]       { $g_{17}$};
  
  \node (20r) [ right       of=20]       {};

  \node[main node] (21) [ right      of=20r]       { $g_{4}^{\mathbb{C}}$};

  \node[main node] (22) [right         of=21]       { $g_{3}^{\mathbb{C}}$};

\node (20b) [below left of=20] {};   

  \node[main node] (23) [ below right        of=20b]       { $n_3 \oplus n_3 $};

 \node (23b) [below of=23]{};
 
   \node[main node] (25) [ right           of=23b]       { $g_{21}$};
   
    \node (25r) [right of=25]{};
        \node (25rr) [right of=25r]{};
 
     \node[main node] (24) [  right          of=25r]       { $n_4 \oplus \mathbb{C}^2 $};

 \node (24b) [below of=25]{};  

  \node[main node] (26) [  left           of=24b]       { $g_{24}$};

  \node[main node] (27) [ below           of=26]       { $g_{2}^{\mathbb{C}}$};

   \node (27b) [below of=27]{};  

  \node[main node] (28) [ right           of=27b]       { $g_{1}^{\mathbb{C}}$};

  \node[main node] (29) [ below           of=27b]       { $n_3 \oplus \mathbb{C}^3$};

  \node[main node] (ab) [ below           of=29]       { $ \mathbb{C}^6$};

  \node (o1) [ left       of=11]       { $12$};
  
  \node (o1a) [above of=o1] {};

  \node (o2) [ above       of=o1a]       { $11$};
  
  \node (o2a) [above of=o2] {};

  \node (o3) [ above       of=o2a]       { $10$};
  
  \node (o3a) [above left of=o3] {};

  \node (o4) [ above right        of=o3a]       { $9$};

  \node (o5) [ above        of=o4]       { $8$};
  
  \node (o1b) [below left of=o1] {};

  \node (o6) [ below  right      of=o1b]       { $13$};
  
  \node (o6b) [below left of=o6] {};

  \node (o7) [ below   right     of=o6b]       { $14$};
  
  \node (o7b) [below left of=o7] {};

  \node (o8) [ below   right     of=o7b]       { $15$};
  
  \node (o8b) [below left of=o8] {};
  
  \node (o9) [ below   right     of=o8b]       { $16$};

  \node (o10) [ below         of=o9]       { $17$};

  \node (o11) [ below          of=o10]       { $18$};

  \node (o12) [ below         of=o11]       { $19$};

  \node (o13) [ below          of=o12]       { $21$};

  \node (o14) [ below        of=o13]       { $24$}; 
  
    \node (o15) [ below        of=o14]       { $36$};

  \path[every node/.style={font=\sffamily\small}]
    (1) edge   node[above] {} (2)

    (1) edge   node[left] {} (3)

    (2) edge   node[left] {} (4)
    
    (3) edge   node[left] {} (4)

    (3) edge   node[left] {} (5)

    (3) edge  node[right, fill=white]{\tiny $\epsilon=1$} (6)
    
    (2) edge  node[above, fill=white]{\tiny $\epsilon=1$}   node[left] {} (9)
    
    (2) edge  [bend left=50] node[left] {} (10)
    
    (4) edge   node[left] {} (7)

    (4) edge   node[left] {} (8)

    (4) edge   node[left] {} (13)

    (4) edge [bend right=20]  node[left] {} (14)

    (5) edge   node[left] {} (8)

    (5) edge   node[above, fill=white]{\tiny $\epsilon=1$} node[left] {} (9)
    
    (5) edge   node[left] {} (10)
    
    (6) edge   [bend right=10] node[left] {} (9)
    
    (6) edge   [bend right=10] node[left] {} (10)

    (8) edge   node[left] {} (11)

    (8) edge   node[left] {} (12)

    (9) edge   node[left] {} (13)

    (10) edge   node[left] {} (11)

    (10) edge   node[left] {} (12)

    (10) edge   node[left] {} (13)
    
    (10) edge   node[left] {} (14)

    (7)  edge   node[left] {} (15)

    (11) edge   node[left] {} (15)

    (11) edge   node[left] {} (16)

    (12) edge   node[left] {} (16)

    (12) edge   node[left] {} (17)
    
    (13) edge   node[left] {} (17)

    (13) edge   node[left] {} (19)

    (14) edge   node[left] {} (17)
    
    (14) edge   node[left] {} (15)

    (16) edge   node[left] {} (18)

    (17) edge   node[left] {} (18)

    (14) edge   node[left] {} (19)

    (15) edge   node[left] {} (20)

    (18) edge   node[left] {} (20)

    (16) edge   node[left] {} (21)

    (14) edge   node[left] {} (21)

    (18) edge   node[left] {} (22)
    
    (19) edge   node[left] {} (22)
    

    (19) edge   node[left] {} (23)

    (17) edge    node[left] {} (23)

    (20) edge   node[left] {} (24)

 
    (21) edge    node[left] {} (24)

    (22) edge  node[left] {} (24)

    (22) edge   node[left] {} (25)

    (23) edge   node[left] {} (25)

    (21) edge     (26)

    (18) edge [bend right] (26)

    (26) edge   node[left] {} (27)

    (25) edge   node[left] {} (27)

    (24) edge    node[left] {} (27)

    (25) edge   node[left] {} (28)

    (28) edge   node[left] {} (29)

    (27) edge   node[left] {} (29)

    (29) edge   node[left] {} (ab)

    (m5c) edge   node[left] {} (22)
    
     (m20) edge   node[left] {} (m5c)

   (m111) edge   node[left] {} (m101)

  (m111) edge   node[left] {} (m110)
  
    (m101) edge   node[left] {} (m5c)
    
      (m110) edge   node[left] {} (m5c)   
   
   (6) edge node[above, rectangle, rounded corners=1.5ex, fill=white]{\tiny $\epsilon=0$} (m50)
   
   (6) edge [bend right=21]  node[above, fill=white]{\tiny $\epsilon=0$} (m110)
   
   (6) edge [bend left=58] node[above, fill=white]{\tiny $\epsilon=0$} (m4)
   
      (6) edge [bend left=4] node[above, rectangle, rounded corners=2.5ex, fill=white]{\tiny $\epsilon=-1$} (m3)
      
      (m71) edge (10)
      
      (m71) edge  (m4)
      
    (m71) edge (m70)
    
    (m71) edge (m111)
        
    (m70) edge (m101)
  
    
    (m70) edge (m20)

    (m51) edge (m50)
    
    (m50) edge (m20)
    
    (m4) edge (m5c)
    
    (m3) edge (m2-1)
    
    (m50) edge (12)
    
    (m20) edge (17)
    
    (m3) edge (13)
    
    (m2-1) edge (18)
    
    (m5c) edge [bend left] (23)
    
    (m70)  edge [bend right=26] (15)
    
    (m101) edge (19)
    
    (m110) edge [bend right=10] (19)
    
    (m4) edge (18)
    
    (m4) edge (19)
    
    (m70) edge [bend left=20] (19);
\end{tikzpicture}

\

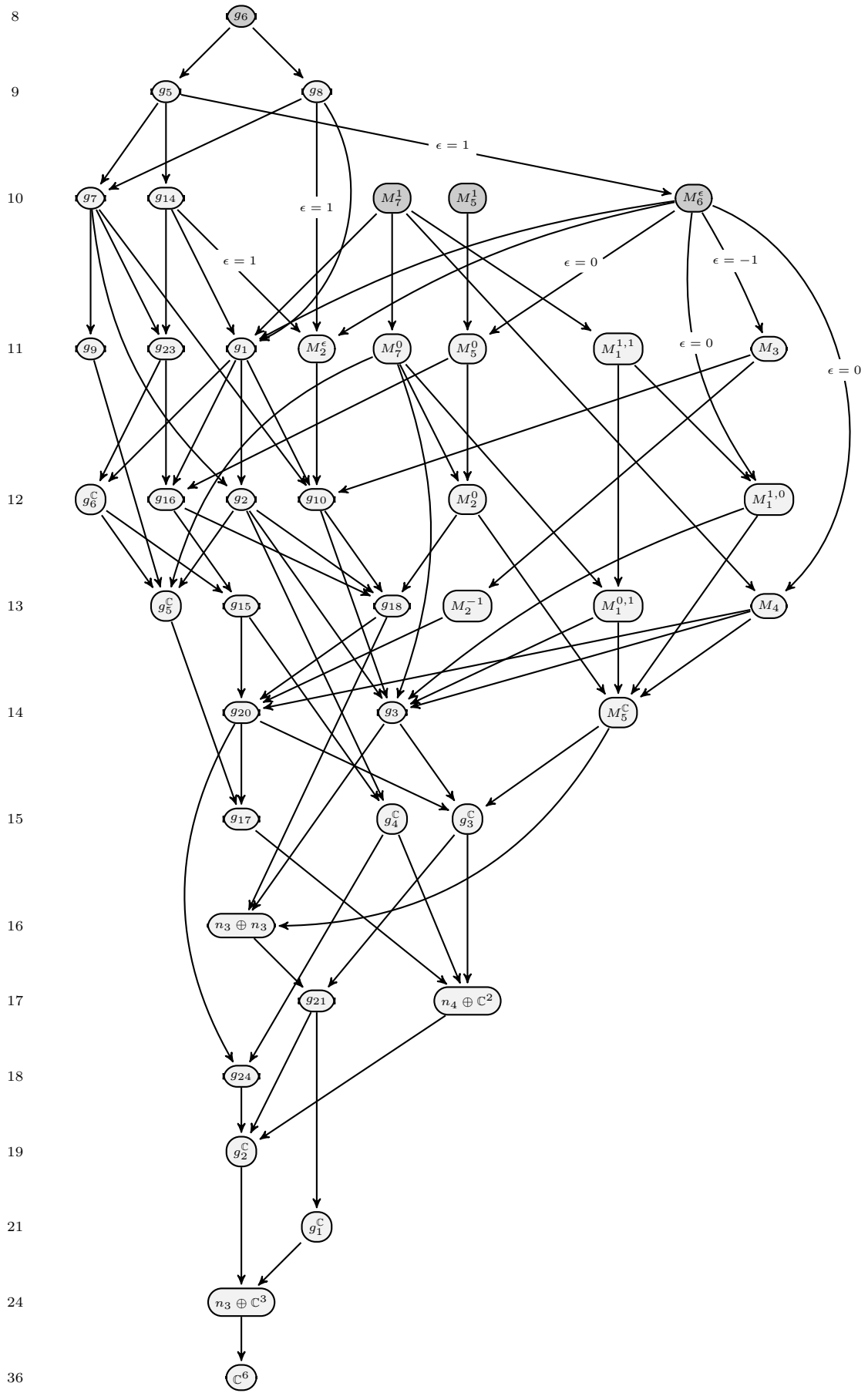
{\bf Figure III.} The graph of primary degenerations for six-dimensional nilpotent Malcev algebras.
\end{center}

\begin{corollary}
$NMal_6=\{\mathcal{C}_1,\mathcal{C}_2\}$, where
$\mathcal{C}_1=\overline{O(g_6)}=NLie_6$ and $\mathcal{C}_2=\overline{\bigcup\limits_{\epsilon\in\mathbb{C}}O(M_6^{\epsilon})}=NMal_6\setminus\{g_6,g_5,g_8,g_7,g_{14},g_9,g_{23}\}$
In particular, $Rig(NMal_6)=Rig(NLie_6)=\{g_6\}$.
\end{corollary}
\begin{Proof} In view of Theorem \ref{third} it is enough to prove that $M_6^*\not\to g_9$, $M_6^*\not\to g_{23}$, $M_6^*\to M_7^1$ and $M_6^*\to M_5^1$, where $M_6^*=\{M_6^{\epsilon}\}_{\epsilon\in\mathbb{C}}$.
The first two assertions follow from the fact that $dim\,(g_9)^2=dim\,(g_{23})^2>dim\,A^2$ for any $A\in M_6^*$.

To prove that $M_6^*\to M_7^1$ one can choose the parametrized basis
$$
E_1^t=e_1,E_2^t=e_2-e_4,E_3^t=te_4,E_4^t=e_3,E_5^t=e_5,E_6^t=e_6
$$
and the parametried index $\epsilon(t)=\frac{1}{t}$.

To prove that $M_6^*\to M_5^1$ one can choose the parametrized basis
$$
E_1^t=e_2,E_2^t=te_1,E_3^t=e_4,E_4^t=-te_3,E_5^t=-t^2e_5,E_6^t=te_6
$$
and the parametrized index $\epsilon(t)=-t^2$.
\end{Proof}

{\bf Acknowledgements.}
The authors are grateful to Prof. Dietrich Burde for some constructive comments.


\begin{thebibliography}{99}

\bibitem{AOR05} 
Albeverio S., Omirov B.,  Rakhimov I., 
Varieties of nilpotent complex Leibniz algebras of dimension less than five, 
Comm. Algebra, 33 (2005), 5, 1575--1585. 





\bibitem{BB09} Benes T., Burde D., Degenerations of pre-Lie algebras, J. Math. Phys., 50 (2009), 11, 112102.

\bibitem{BB14} Benes T., Burde D., Classification of orbit closures in the variety of three-dimensional Novikov algebras, J. Alg. Appl., 13 (2014), 2, 1350081.



\bibitem{B99} Burde D., Degenerations of nilpotent Lie algebras,
J. Lie Theory, 9 (1999), 1, 193--202.

\bibitem{B03} Burde D.,
Sur les degenerations d'algebres de Lie,
\href{https://homepage.univie.ac.at/Dietrich.Burde/papers/burde_15_rapp_deg.pdf}{https://homepage.univie.ac.at/Dietrich.Burde/papers/burde\_15\_rapp\_deg.pdf} (2003).



\bibitem{B05} Burde D., Degenerations of 7-dimensional nilpotent Lie algebras,
Comm. Algebra, 33 (2005), 4, 1259--1277.

\bibitem{BC99} Burde D., Steinhoff C., Classification of orbit closures of 4--dimensional complex Lie algebras,  
J. Algebra, 214 (1999), 2, 729--739.


\bibitem{CKLO13}
Casas J., Khudoyberdiyev A., Ladra M., Omirov B., On the degenerations of solvable Leibniz algebras, 
Lin. Alg. Appl.,  439 (2013),  2, 472--487



\bibitem{G75} 
Gabriel P., Finite representation type is open, 
in: Proceedings of the International Conference on Representations of Algebras, Carleton University, Ottawa, Ontario, 1974, in: Lecture Notes in Math., vol. 488, 1975, pp. 132--155. 

\bibitem{G63}
 Gainov A., Binary Lie algebras of lower ranks (Russian), Algebra i Logika Sem., 2 (1963), 4, 21--40.

\bibitem{gai57}
Gainov A., Identical relations for binary Lie rings (Russian), Uspehi Mat. Nauk N.S., 12 (1957), 3 (75), 141--146. 




\bibitem{gorb93}
Gorbatsevich V., Anticommutative finite-dimensional algebras of the first three levels of complexity,
St. Petersburg Math. J., 5 (1994), 3, 505--521.	



\bibitem{GRH}
Grunewald F.,  OВґHalloran J., 
Varieties of nilpotent Lie algebras of dimension less than six, J. Algebra, 112 (1988), 315--325.

\bibitem{GRH2}
Grunewald F., OВґHalloran J., 
A Characterization of Orbit Closure and Applications, J. Algebra, 116 (1988), 163--175.




\bibitem{fil2}
Filippov V., On $\delta$-derivations of prime alternative and Malcev algebras,  
Algebra and Logic, 39 (2000), 5, 618--625.


\bibitem{KE14}  Kashuba I., Martin M., Deformations of Jordan algebras of dimension four, J. Algebra, 399 (2014), 277--289.


\bibitem{KE16}  Kashuba I., Martin M., The variety of three-dimensional real Jordan algebras,
J. Alg. Appl.,  15 (2016), 8, 1650158.


\bibitem{kay14}
Kaygorodov I., On $(n+1)$-ary derivations of simple $n$-ary Malcev algebras, 
St. Petersburg Math. J., 25 (2014), 4, 575--585

\bibitem{kp16}
Kaygorodov I., Popov Yu.,
A characterization of nilpotent nonassociative algebras by invertible Leibniz-derivations,
J. Algebra, 456 (2016), 323--347.


\bibitem{K70}
Kuzmin E., 
Malcev algebras of dimension five over a field of characteristic zero,
Algebra and Logic, 9 (1970), 416--421.

\bibitem{kuz71}
Kuzmin E., The connection between Malcev algebras and analytic Moufang loops, 
Algebra and Logic, 10 (1971), 3--22.



\bibitem{K98}
Kuzmin E., Binary Lie algebras of small dimension, Algebra and Logic, 37 (1998), 3, 181--186.


\bibitem{kuzma2}
Kuzmin E., Structure and representations of finite dimensional Malcev algebras, 
Quasigroups Related Systems, 22 (2014), 1, 97--132.




\bibitem{M55}
Malcev A., Analytic loops (Russian), Mat. Sb. N.S., 36, (1955), 569--576.

\bibitem{M79} 
Mazzola G., The algebraic and geometric classification of associative algebras of dimension five, 
Manuscripta Math., 27 (1979), 81--101. 

\bibitem{M80}  
Mazzola G., Generic finite schemes and Hochschild cocycles, Comment. Math. Helv. 55 (1980), 267--293. 

\bibitem{poji01}
Pozhidaev A., $n$-ary Malcev algebras,  Algebra and Logic, 40 (2001),  3, 170--182.

\bibitem{R06} 
Rakhimov I., On the degenerations of finite dimensional nilpotent complex Leibniz algebras, 
J. Math. Sci. (N.Y.), 136 (2006), 3, 3980--3983.


\bibitem{sagle61}
Sagle A., Malcev algebras, Trans. Amer. Math. Soc., 101 (1961), 426--458.


\bibitem{S90}
Seeley C., Degenerations of 6-dimensional nilpotent Lie algebras over $\mathbb{C}$, Comm. Algebra, 18 (1990), 3493--3505.




\end{thebibliography}
\end{document}